\documentclass{amsart}

\usepackage{amsmath,amssymb,amsthm,graphicx,slashed,subfig}%
\usepackage{amsmath}%
\usepackage{amsfonts}%
\usepackage{amssymb}%

\usepackage{mathpple}

\usepackage{hyperref,enumerate}

\usepackage{colortbl}

\usepackage[all,cmtip]{xy}

\definecolor{Red}{rgb}{0.6,0,0}

\usepackage{graphicx}
\providecommand{\U}[1]{\protect\rule{.1in}{.1in}}
\newtheorem{thm}{Theorem}[section]
\newtheorem{corl}[thm]{Corollary}

\newtheorem{lma}[thm]{Lemma}

\newtheorem{prop}[thm]{Proposition}
\newtheorem{defn}[thm]{Definition}
\newtheorem{ex}[thm]{Example}
\newtheorem{rem}[thm]{Remark}
\def\tilde{\widetilde}

\def\A{\mathcal{A}}
\def\cS{\mathcal{S}}
\def\B{\mathcal{B}}
\def\bar{\overline}
\newcommand{\bra}[1]{\langle #1 |}
\newcommand{\ket}[1]{| #1 \rangle}
\def\C{\mathbb{C}}
\def\cb{\mathrm{cb}}
\def\CB{\mathrm{CB}}

\def\E{\mathcal{E}}
\def\env{\mathrm{env}}
\def\epsilon{\varepsilon}
\def\eps{\varepsilon}

\renewcommand{\Im}{\textup{Im}}

\def\id{\mathbb{I}}
\def\cK{\mathcal{K}}
\def\L{\mathcal{L}}

\def\min{\mathrm{min}}

\def\cP{\mathcal P}
\def\P{\mathbb P}

\def\phi{\varphi}

\renewcommand{\Re}{\textup{Re}}
\def\cR{\mathcal{R}}

\DeclareMathOperator{\supp}{\textup{supp}}
\def\cS{\mathcal{S}}

\DeclareMathOperator{\tr}{Tr}
\newcommand{\Toep}[1]{C(S^1)^{(#1)}}
\def\U{\mathcal{U}}
\def\Z{\mathbb{Z}}

\newtheoremstyle{commentstyle}
  {0.2cm}{0.2cm}
  {\sf}
  {0cm}
  {\bfseries}{ }
  {0cm}
  {\thmname{#1}\thmnumber{ #2}:\thmnote{ #3}}

\theoremstyle{commentstyle}
\newtheorem{mycomment}{Comment}

\title{Tolerance relations and operator systems}
\author{Alain Connes}
\address{College de France, 3 rue Ulm, F75005, Paris, France\\
I.H.E.S. F-91440 Bures-sur-Yvette, France\\
}
\email{alain@connes.org}

\author{Walter D. van Suijlekom}

\address{Institute for Mathematics, Astrophysics and Particle Physics, Radboud
University Nijmegen, Heyendaalseweg 135, 6525 AJ Nijmegen, The Netherlands.
}
\email{waltervs@math.ru.nl}
\date{3 November 2021}

\begin{document}

\maketitle
\begin{abstract}
  We extend the scope of noncommutative geometry by generalizing the construction of the noncommutative algebra of a quotient space to situations in which one is no longer dealing with an equivalence relation. For these so-called tolerance relations, passing to the associated equivalence relation looses crucial information as is clear from the examples such as coarse graining in physics or the relation $d(x,y)< \epsilon$ on a metric space. Fortunately, thanks to the formalism of operator systems such an extension is possible and provides new invariants, such as the $C^*$-envelope and the propagation number.

  After a thorough investigation of the structure of the (non-unital) operator systems associated to tolerance relations, we analyze the corresponding state spaces. In particular, we determine the pure state space associated to the operator system for the relation $d(x,y)< \epsilon$ on a path metric measure space.

 \end{abstract}

 \tableofcontents

\section{Introduction}
\label{sect:intro}

In this paper we take the next step in the development of noncommutative geometry based on operator systems. While in \cite{CS20} we considered operator systems based on spectral truncations, we will now focus on operator systems associated to tolerance relations.

A strong motivation for studying such systems comes from the idea of coarse graining in physics, which may be formulated in phase space or in configuration space. Moreover, 
in our understanding of the geometry of spacetime we are limited by the resolution of our measuring device: we can only determine the underlying metric space {\em up to a finite resolution}. In fact, on purely theoretical grounds one expects our understanding of this geometry to break down at a fundamental scale, given by the Planck length $\ell_P = \sqrt{{\hbar G}/{c^3}}$.

From a mathematical viewpoint this leads us to consider metric spaces $(X,d)$ equipped with the relation $d(x,y) <\eps$ for some fixed parameter $\eps$. This is not an equivalence relation as it lacks transitivity, but it is a so-called {\em  tolerance relation}, which is by definition a reflexive and symmetric relation. This notion can be traced back to Poincar\'e in his {\em Science and Hypothesis} (though the name was coined in \cite{Zee62}). In fact, Poincar\'e argued that in the physical continuum (in contrast with the mathematical continuum) it can hold that for measured quantities one has $A=B$, $B=C$ while $A \neq C$ due to potentially added measurement errors (see \cite{Sos86} for a development of the mathematical theory).

Given such a tolerance relation we may imitate the construction of a groupoid $C^*$-algebra with the crucial difference that the convolution product cannot be defined due to the lack of transitivity. However, given the symmetry and reflexivity, a tolerance relation does  define an operator system. Naturally, this operator system may be formulated in the general context of locally compact groupoids equipped with a Haar system \cite{Ren80} with the special case of the tolerance relation arising when we consider the pair groupoid on a locally compact topological measure space. It creates a generalization of the basic construction of noncommutative geometry which started from analyzing the geometric examples of intractable spaces of leaves of foliations using noncommutative algebras. Here the issue of the lack of transitivity of the relation is already present in the simplest case where the generated equivalence relation has a single class {\em i.e.} corresponds to the  $C^*$-algebra of compact operators. The main question is then how much of the metric structure of $(X,d)$ can be determined from the operator system associated to the relation $d(x,y)<\eps$.\newline
Another natural appearance of tolerance relations occurs in combinatorial topology. The homotopy groups of the geometric realization of a Kan complex may be computed combinatorially. This no longer holds for general simplicial complexes and one motivation for considering tolerance relations comes from the natural examples of ``sets" corresponding to a finer notion of the set of components of a simplicial complex $X$ which is not a Kan complex. Typically the relation among  two vertices $x,y\in X(0)$ given by the existence of an edge $e\in X(1)$ whose end points are $x,y$ is a tolerance relation which only uses the graph underlying the simplicial complex in degrees $0,1$. On the one hand replacing the relation by  the equivalence relation that it generates and then passing to the quotient looses most of the relevant information. On the other hand one would like to perform this operation of quotient for factor-relations which are equivalence relations. It is this issue which is solved by the construction of the associated operator system and the notion of Morita equivalence of operator systems which was developed recently in \cite{EKT21}.

A concrete example of this situation occurs in the theory of divisors on the Arakelov compactification $\overline{{\rm Spec} \; \mathbb Z}$ (see \cite{CC21}). After applying the Dold--Kan correspondence to a suitable morphism of $S$-modules, one obtains that the cohomology $H^1(\overline{{\rm Spec} \;\mathbb Z},D)$ should be described by the tolerance relation given by the circle $S^1$ endowed with the relation 
$$
\cR:=\{(x,y)\mid d(x,y)<\epsilon\}
$$
where $\epsilon=\exp(\deg D)$ becomes small when the degree $\deg D$ of the divisor is negative.

\medskip

This paper is organized as follows. We start in Section \ref{sect:opsyst} by reviewing the theory of non-unital operator systems, due to Werner \cite{Wer02}, as well as the notions of $C^*$-envelopes and propagation number as developed previously in \cite{CS20}. Under the assumption that the operator system $E$ is a closed $*$-subspace of a $C^*$-algebra $A$ equipped with a suitable approximate order unit (in the sense of \cite{Ng69}), we prove a generalization of the Jordan decomposition for hermitian functionals on $E$ (Theorem \ref{thm:jordan}). This then allows us to characterize the Banach space dual of $E$ in terms of the dual of $A$, preparing for an analysis of the corresponding state space. 

In Section \ref{sect:opsyst-bonds} we define operator systems associated to relations. We work in the slightly more general, but natural context of locally compact groupoids equipped with a Haar system, and associate operator systems to open symmetric subsets (so-called {\em bonds}, see Definition \ref{defn:bonds}). After discussing some of their general properties, we analyze operator systems associated to tolerance relations on finite sets, and determine their $C^*$-envelopes, propagation numbers and pure state space. 

Finally, in Section \ref{sect:relations} we come to consider the operator system associated to the relation $d(x,y) < \eps$ on a path metric space $X$. In addition to invariants such as the $C^*$-envelope and the propagation number (Theorem \ref{thm:prop-nr}), we determine the pure state space of this operator system to consist of those vector states whose support is $\eps$-connected (Theorem \ref{prop:pure2}).

\section{Preliminaries on non-unital operator system and state spaces}
\label{sect:opsyst}

We start by recalling and developing some general theory for non-unital operator systems, approximate order units and state spaces.


\subsection{Non-unital operator systems}
Werner introduced non-unital operator system in \cite{Wer02}. Let us give the relevant definitions, referring to {\em loc.cit.} and also \cite{CS20} for more details and references. 

\begin{defn}\label{defn:matrix-order}
We say that a $*$-vector space $E$ is {\em matrix ordered} if
\begin{enumerate}
\item for each $n$ we are given a cone of positive elements $M_n(E)_+$ in $M_n(E)_h$,
\item $M_n(E)_+ \cap (-M_n(E)_+) = \{ 0\}$ for all $n$,
  \item for every $m,n$ and $A  \in M_{mn}(\C)$ we have that $A M_n(E)_+ A^* \subseteq M_m(E)_+$.
\end{enumerate}
\end{defn}

For {\em matrix-ordered operator spaces} one further requires that $E$ is an operator space, with the two properties that
\begin{enumerate}
\item the cones $M_n(E)_+$ are all closed, and,
\item the involution is an isometry on $M_n(E)$.
\end{enumerate}

The {\em modified numerical radius} $\nu_E$ is defined for any $x \in M_n(E)$ by 
$$
\nu_E(x)  = \sup \left \{ \left | \phi \begin{pmatrix} 0 & x \\ x^* & 0 \end{pmatrix} \right| : \phi \in M_{2n}(E)^*_+, \| \phi \| \leq 1 \right \}. 
$$
Note that $M_{2n}(E)^*$ is the Banach space dual of $M_{2n}(E)$, in contrast to $M_{2n}(E^*) = \CB(E, M_{2n})$ which comes naturally equipped with the $\cb$-norm. The supremum is thus taken over all noncommutative quasi-states of $E$, {\em cf.} Equation \eqref{eq:nc-quasi-states} below.

One can show that $\nu_E$ is a norm and, moreover, we have $\nu_E(\cdot) \leq \| \cdot \|$. Werner then gives the following definition in \cite{Wer02} (where we prefer as in \cite{CS20} to restrict to identical, rather than equivalent norms):
\begin{defn}
  A {\em non-unital operator system} is given by a matrix-ordered operator space $E$ for which the norm $\nu_E(\cdot)$ coincides with the norm $\| \cdot \|$.
\end{defn}

Werner then constructs a unitization for arbitrary matrix-ordered operator spaces as follows.
\begin{defn}
\label{defn:unit-opsyst}
Let $E$ be a matrix-ordered operator space and define $a_\eps = a + \eps \id_n$ for every matrix $a \in M_n(\C)$. On the space $E \oplus \C$ we define
\begin{enumerate}
\item $(x,a)^* = (x^*,a^*)$ for all $(x,a) \in M_n(E )\oplus M_n( \C)$,
\item for any $(x,a) \in M_n(E \oplus \C)_h$ we set
$$ (x,a) \geq 0 \qquad \text{iff } \quad a \geq 0 \text{ and } \phi(a_\eps^{-1/2} x a_\eps^{-1/2}) \geq -1
$$
for all $\eps>0$ and $\phi \in M_n(E)^*_+$ with $\| \phi\|=1$.
\end{enumerate}
We denote by $E^\sharp$ the space $E \oplus \C$ equipped with this order structure. 
\end{defn}

The main conclusion from \cite{Wer02} is then that $E^\sharp$ is a unital operator system and that $E$ is a non-unital operator system iff the embedding $\imath : E \to E^\sharp$ is a complete isometry. For a convex geometric description of operator systems and their unitization we refer to \cite{DK15,KKM21}. The following result is a special case of \cite[Proposition 4.16(a)]{Wer02}
\begin{lma}
  \label{lma:opsyst-Cstar}
  Let $E \subseteq A$ be a closed $*$-subspace in a $C^*$-algebra $A$. Then $E$ is a non-unital operator system.
  \end{lma}
\proof
Take $\left(\begin{smallmatrix} 0 & x \\ x^* & 0 \end{smallmatrix} \right) \in M_{2n}(E)$; this is a self-adjoint element in $M_{2n}(A)$. Hence there exists a state $\tilde \phi$ on the $C^*$-algebra $M_{2n}(A)$ such that $\left|\tilde \phi \left(\left(\begin{smallmatrix} 0 & x \\ x^* & 0 \end{smallmatrix} \right) \right) \right| =  \| x \|$. The restriction $\phi$ of this state to $M_{2n}(E)$ is a quasi-state and hence $\nu_E(x) \geq |\phi(x) | = \| x\|$. 
\endproof

In \cite{CS20} we introduced the notion of {\em propagation number} of an operator system, a measure for the extent to which the operator system is `away' from a $C^*$-algebra. Let us start by briefly recalling the notion of the $C^*$-envelope \cite{Arv69,Ham79}

\begin{defn}
  \begin{enumerate}
  \item
  A {\em $C^*$-extension} $\kappa: E \to A$ of a unital operator system $E$ is given by a complete order isomorphism onto $\kappa(E) \subseteq A$ such that $C^*(\kappa(E)) = A$.
  \item A {\em $C^*$-envelope} of a unital operator system is a $C^*$-extension $\kappa: E \to A$ with the following universal property: 
    $$
    \xymatrix{
      E \ar^\kappa[r]\ar_\lambda[rd] &A \ar^{\exists ! \rho}@{<<.}[d]\\
      & B 
      }$$
\end{enumerate}
\end{defn}
Such a universal object was shown to exist in the unital case in \cite{Ham79} while for the non-unital case we refer to \cite{CS20} and references therein. We write $\imath_E : E \to C^*_\env(E)$ for the canonical complete order injection into the $C^*$-envelope. 
      
Given an operator system $E$ and integer $n>0$ let us write $E^{\circ n}$ for the norm closure of the linear span of products of $\leq n$ elements of $E$. It is an operator system contained in the $C^*$-algebra $C^*(E)$. 
\begin{defn}
The {\em propagation number} ${\rm prop}(E)$ of the operator system $E$ is defined as the smallest integer $n$ such that $\imath_E(E)^{\circ n}\subseteq C_{\rm env}^*(E)$ is a $C^*$-algebra. 
  \end{defn}
When no such $n$ exists one lets ${\rm prop}(E)=\infty$. In \cite[Proposition 2.42]{CS20} we showed that the propagation number is invariant under stable equivalence of operator systems. Recently, also a notion of Morita equivalence for (unital) operator systems has been introduced \cite{EKT21} and shown to coincide with stable equivalence (see Theorem 3.8 in {\em loc.cit.}). We may thus conclude that the propagation number is also invariant under Morita equivalence. Other properties of the propagation number include compatibility with tensor products. Namely, in \cite{Koo21} it is shown that for unital operator systems $E$ and $F$ we have ${\rm prop}(E \otimes_{\text{min}} F ) = \max \{{\rm prop}(E), {\rm prop}(F)\}$.

\subsection{Matrix-regular operator spaces and approximate order units}
We now introduce a certain degree of regularity on matrix-ordered operator spaces, which relates the matrix norms and the matrix order, in fact yielding a special class of non-unital operator systems as considered above. A key concept will be that of an approximate order unit due to Ng \cite{Ng69}, but let us start with the following, more general notion introduced by Schreiner in \cite{Sch98}. 

\begin{defn}
  A matrix-ordered operator space $E$ is called a {\em matrix-regular operator space} if for all $x \in M_n(E)_h$ the following two conditions are satisfied
  \begin{enumerate}
  \item if $y \in M_n(E)_h, -y \leq x \leq y$ then $\|x \| \leq \| y \|$;
    \item if $\| x \| <1$ then there exists $y \in M_n(E)_h$ such that $\| y \| \leq 1$ and $-y \leq x \leq y$. 
    \end{enumerate}
\end{defn}
Note the equivalent characterization of matrix-regular operator spaces \cite[Theorem 3.4]{Sch98}: for all $x \in M_n(E)$, $\| x \| < 1$ if and only if there exist $a,d\in M_n(E)_+$, $\| a \| , \| d \|<1$ such that
$$
\begin{pmatrix} a & x \\ x^* & d \end{pmatrix} \in M_{2n}(E)_+.
$$
Schreiner establishes in \cite[Corollary 4.7]{Sch98} that the operator space dual of a matrix-regular operator space is also matrix regular. Moreover, Karn shows in \cite{Kar05,Kar07} ({\em cf.} \cite[Theorem 1]{Han10}) that matrix-regular operator spaces are non-unital operator systems. In this context we should also mention the recent work on dual spaces of non-unital operator systems in \cite{Ng21}.

\medskip 
A particular class of matrix-regular operator spaces is given by operator spaces that possess an approximate order unit in the sense of Ng in \cite{Ng69}.
\begin{defn}
  \label{defn:approx-order-opsyst}
Let $E$ be a matrix-ordered $*$-vector space. An {\em approximate order unit} for $E$ is a net $\{e_\lambda \}_{\lambda \in \Lambda}$ with the following properties
\begin{enumerate}
\item $e_\lambda \in E_+$ for all $\lambda \in \Lambda$;
\item $e_\lambda \leq e_\mu $ whenever $\lambda \leq \mu$;
\item for each $x \in E_h$ there exists a positive real number $t$ and $\lambda \in \Lambda$ such that
  $$
  - t e_\lambda \leq x \leq t e_\lambda.
  $$
 \end{enumerate}
\end{defn}


\begin{lma}{\cite[Lemma 2.6]{KV97}}
Let $E$ be a matrix-ordered $*$-vector space. If $E$ has an approximate order unit, then $M_n(E)$ has an approximate order unit for all $n$.
\end{lma}
\proof
Our proof is inspired by the approach taken in the unital case \cite[Theorem 4.4]{CE77}. First of all, for $x \in E_h$ we have ({\em cf.} \cite[Proposition 1.8(i)]{KV97}
$$
-te_\lambda \leq x \leq te_\lambda \iff \begin{pmatrix} t e_\lambda & x \\ x & t e_\lambda \end{pmatrix} \geq 0.
$$
For a general $x \in E$ we may apply this to both $\Re(x), \Im(x)$ to conclude that
$$
\begin{pmatrix} 2t e_\lambda & x \\ x^* & 2t e_\lambda \end{pmatrix} = 
\begin{pmatrix} t e_\lambda & \Re(x) \\ \Re(x) & t e_\lambda \end{pmatrix}
+ 
\begin{pmatrix} 1 & 0 \\ 0 & -i \end{pmatrix} \begin{pmatrix} t e_\lambda & \Im(x) \\ \Im(x) & t e_\lambda \end{pmatrix}
\begin{pmatrix} 1 & 0 \\ 0 & i \end{pmatrix}\geq 0.
$$
Thus for $x = (x_{ij}) \in M_n(E)_h$ we may choose $t$ and $\lambda\in \Lambda$ such that $-t e_\lambda \leq \Re(v_{ij}), \Im(v_{ij} \leq t e_\lambda$ for all $i,j$. Hence
$$
2 n t  e_\lambda^n + x = \frac 12 \sum_{i,j} A_{ij}^*  \begin{pmatrix} 2t e_\lambda & x_{ij} \\ x_{ji} & 2t e_\lambda \end{pmatrix} A_{ij} \geq 0,
$$
where $A_{ij}$ is the $2 \times m$ matrix with $1$ at the $(1,i)$'th and $(2,j)$'th positions and 0 elsewhere for $i,j=1,\ldots, n$. Similarly, we find that $2 n t  e_\lambda^n - x \geq 0$ which proves the claim.
\endproof

As we see from the proof we have that for all $x \in M_n(E)$ there exist $t>0$ and $\lambda \in \Lambda$ such that 
\begin{equation}
    \label{eq:pos-approx}
\begin{pmatrix} t e_\lambda^n & x \\ x^* & t e_\lambda^n \end{pmatrix} \in M_{2n}(E)_+
\end{equation}
This is the key property that allows us to introduce a matrix norm on $E$ (as in \cite{Ng69}, \cite{Kar05,Kar07}). Indeed, given an approximate order unit $\{ e_\lambda\} _{\lambda \in \Lambda}$ we may introduce a seminorm $p_{\Lambda,n}$ on $M_n(E)$ by setting 
  $$
p_{\Lambda,n}(x) = \inf \left\{ t :  \begin{pmatrix} t e_\lambda^n & x \\ x^* & t e_\lambda^n \end{pmatrix}  \in M_{2n}(E)_+ \text{ for some } \lambda \in \Lambda \right\}
$$
In the case of a matrix-ordered operator space, which is equipped with an approximate order unit, it is natural to require the seminorms $p_{\Lambda,n}$ to be compatible with the matrix norms, which leads to the following result.

\begin{prop}
  \label{prop:approx-ord-matrix-regular}
Suppose that $E$ is a matrix-ordered operator space with an approximate order unit that defines the matrix norms (in the sense that $p_{\Lambda,n} = \| \cdot \|$). Then $E$ is matrix regular. Consequently, $E$ is a non-unital operator system.
  \end{prop}
\proof 
Let $x \in M_n(E)$ with $\| x \|<1$. Then there exist $t<1$ so that
$$
\begin{pmatrix} t e_\lambda^n & x \\ x^* & t e_\lambda^n \end{pmatrix} \in M_{2n}(E)_+.
$$
where indeed $t e_\lambda^n \geq 0$ and $\| t e_\lambda^n \| < 1$.

Conversely, suppose there exist $a,d\in M_n(E)_+$, $\| a \| , \| d \|<1$ for which $\left(\begin{smallmatrix} a & x \\ x^* & d \end{smallmatrix} \right)\geq 0$. Then also $\left(\begin{smallmatrix} a & -x \\ -x^* & d \end{smallmatrix}\right) \geq 0$ so that
$$
- \begin{pmatrix} a & 0 \\ 0 & d \end{pmatrix}\leq \begin{pmatrix}0  & x \\ x^* & 0 \end{pmatrix} \leq \begin{pmatrix} a & 0 \\ 0 & d \end{pmatrix}.
$$
Since $p_{\Lambda,2n} \left (\left(\begin{smallmatrix} a & 0 \\ 0 & d \end{smallmatrix}\right) \right) < 1$ we thus find for some $t<1$ that
$$
-  t e_\lambda^{2n} \leq  \begin{pmatrix} a & 0 \\ 0 & d \end{pmatrix}  \leq t e_\lambda^{2n}, \implies - t e_\lambda^{2n} \leq  \begin{pmatrix} 0 & x \\ x^* & 0 \end{pmatrix}  \leq t e_\lambda^{2n},
$$
and thus $\| x \| < 1$ as required.

The fact that then $E$ is a non-unital operator systems follows from \cite{Kar05,Kar07} ({\em cf.} \cite[Theorem 1]{Han10}).
\endproof
Note that this result also applies without the assumption of $E$ being complete.

\subsection{States, quasi-states and pure states of matrix-ordered operator spaces}
\label{sect:quasi-states}
Let $E$ be a matrix-ordered operator space. We define the noncommutative (nc) state space for any $n$ as follows:
$$
\cS_n(E) := \{ \phi \in M_n(E)^*, \| \phi \| = 1, \phi \geq 0\}
$$
Note that we consider $M_n(E)^*$ as the Banach space dual of the normed space $M_n(E)$, and not as the operator space dual $M_n(E^*) = \CB(E,M_n)$ which is equipped with the $\cb$-norm (in contrast to e.g. \cite[Section 5.1]{ER00}). 

\begin{lma}
  \label{lma:norm-limit}
  Let $E$ be a matrix-ordered operator space and $\E \subseteq E$ a dense $*$-subspace. Suppose there exists an approximate order unit $\{ e_\lambda \}_{\lambda \in \Lambda}$ in $\E$ that defines the norm on $\E$ and let $\phi \in M_n(E)^*$ for any $n$. If $\phi \geq 0$ then $\| \phi \| = \lim  \phi(e_\lambda^n)$. As a consequence, $\cS_n(E)$ is a convex set. 
\end{lma}
\proof
Clearly, $\|\phi \| \geq \phi(e_\lambda^n)$ for all $\lambda$ so let us prove the other inequality. Take $x \in M_n(\E)$ with $\| x \| \leq 1$. Then we have
$$
- \begin{pmatrix} e_\lambda^n & 0 \\ 0 & e_\lambda^n \end{pmatrix} \leq \begin{pmatrix} 0& x \\ x^* & 0 \end{pmatrix} \leq \begin{pmatrix} e_\lambda^n & 0 \\ 0 & e_\lambda^n \end{pmatrix} .
$$
Since $\E$ is matrix-regular and $\phi \geq 0$ it follows that 
$$
\left \| \begin{pmatrix} 0& \phi(x) \\ \phi(x)^* & 0 \end{pmatrix} \right \| \leq 
\left \| \begin{pmatrix} \phi(e_\lambda^n) & 0 \\ 0 & \phi(e_\lambda^n) \end{pmatrix} \right \| \implies |\phi(x)| \leq   \phi(e_\lambda^n) \leq \liminf \phi(e_\lambda^n).
$$
Taking the supremum over all $x$ with $\| x \| \leq 1$ we find that $\| \phi \| \leq \liminf \phi(e_\lambda^n)$.
\endproof
Note, however, that in general $S_n(E)$ is not weakly $*$-compact. Under the conditions stated in Lemma \ref{lma:norm-limit} we will call an extreme point of $\cS_n(E)$ a {\em (nc) pure state}, and write $\cP_n(E)$ for all (nc) pure states on $E$.

The {\em nc quasi-state space} is defined for any $n$ as
\begin{equation}
  \label{eq:nc-quasi-states}
\tilde \cS_n(E) := \{ \phi \in M_n(E)^*, \| \phi \| \leq 1, \phi \geq 0\}
\end{equation}
In contrast to $S_n(E)$, this is a convex weakly $*$-compact space and contains the point $0$ as an extreme point. Any other extreme point $\phi$ has $\| \phi\|=1$,  for if $\| \phi \| =\lambda <1$ then $\lambda^{-1} \phi \in \tilde \cS(E)$. We record the following result.
\begin{lma}
Let $E$ be a matrix-ordered operator space and $\E \subseteq E$ a dense $*$-subspace containing a norm-defining approximate order unit.   
 Then a state $\phi$ on a matrix-ordered operator space $E$ is pure if and only if it is extreme in $\tilde \cS(E)$
\end{lma}
\proof
If $\phi$ is extreme in $\tilde \cS(E)$ then it cannot be written as a convex combination of elements in $\tilde \cS(E)$, let alone of elements in $\cS(E)$. Conversely, suppose that $\phi$ is not extreme and can thus be written as a convex combination $t\phi_1 + (1-t) \phi_2$ with $\phi_1, \phi_2 \in \tilde \cS(E)$, we claim that in fact $\| \phi_1 \| = \| \phi_2 \| =1$ so that $\phi$ is not pure. Indeed, if on the contrary say $\| \phi_1 \| < 1$ we also have that $\| \phi \| <1$, contradicting the fact that $\phi$ is a state. 
\endproof

\subsection{Extension of positive functionals, Jordan decomposition}
\label{sect:jordan}

Let us start by formulating an Arveson's extension theorem for positive functionals on matrix-ordered operator spaces contained in a $C^*$-algebra, in the presence of a norm-defining approximate order unit. We also refer to \cite[Theorem 3.4]{Kar05} and \cite{Kar07} for an analogous result in the context of matricial Riesz normed spaces.

\begin{prop}
  \label{prop:pos-HB}
  Let $E$ be a closed $*$-subspace in a $C^*$-algebra $A$ and let $\E \subseteq \A$ be dense $*$-subspaces, so that $\E \subseteq E$ and $\A \subseteq A$. If there exists a norm-defining approximate order unit for $\A$ which is contained in $\E$, then any positive linear functional $\phi$ on $M_n(E)$ can be extended to a positive linear functional $\tilde \phi$ to $M_n(A)$ such that $\| \tilde \phi \|  = \| \phi \|$.
    \end{prop}
\proof
By Lemma \ref{lma:opsyst-Cstar} $E$ is a non-unital operator system. Hence, the unitization $E^\sharp$ from Definition \ref{defn:unit-opsyst} is a unital operator system by \cite[Lemma 4.8]{Wer02}.

Now, let $\phi :E \to \C$ be positive and assume that $\| \phi \| =1$. We first consider the extension of $\phi$ to a linear functional $\phi^\sharp : E^\sharp \to \C$ defined by $\phi^\sharp(x \oplus z)= \phi(x) + z$. We claim that the map $\phi^\sharp$ is positive. Indeed, since $(x,z) \geq 0$ iff $z \geq 0$ and $\psi(x) \geq -z$ for all $\psi \in \cS(E)$ it follows that in particular $\phi(x) \geq -z$.

We now apply Arveson's extension Theorem \cite[Theorem 7.5]{Pau02} to the functional $\phi^\sharp$ on $E^\sharp \subseteq A^\sharp$. This gives a positive $\tilde \phi^\sharp: A^\sharp  \to \C$ such that $\tilde \phi^\sharp|_{E^\sharp} = \phi^\sharp$, and thus also $\tilde \phi^\sharp|_{E} = \phi$. We now define $\tilde \phi : A \to \C$ as the restriction $\tilde \phi = \tilde \phi^\sharp|_A$. This map $\tilde \phi$ is the sought-for positive extension of $\phi$ from $E$ to $A$.

Finally, the claim that $\| \tilde \phi \|=  \| \phi\|$ follows from Lemma \ref{lma:norm-limit}. Indeed, we may identify $\E^* \cong E^*$ and $\A^* = A^*$ to find that $\| \tilde \phi \| = \lim \tilde \phi(e_\lambda) = \lim \phi(e_\lambda) = \| \phi \|$ because $e_\lambda \in \E$ is an approximate order unit $\A$.

Note that the generalization to $M_n(E)$ is straightforward since one may consider $M_n(E)$ as a matrix-ordered operator space contained in the $C^*$-algebra $M_n(A)$ and with approximate order unit $e_\lambda^n$.
\endproof
We then have the following well-known result ({\em cf.} \cite[Theorem 5.1.13]{Mur90} adapted to our context as follows.
\begin{prop}
  \label{prop:ext-pure}
Let $E$ be a closed $*$-subspace in a $C^*$-algebra $A$ and let $\E \subseteq \A$ be dense $*$-subspaces, so that $\E \subseteq E$ and $\A \subseteq A$. If there exists a norm-defining approximate order unit for $\A$ which is contained in $\E$, then any pure state $\phi$ on $E$ can be extended to a pure state $\tilde \phi$ on $A$. 
  \end{prop}
\proof
Let $\phi$ be a pure state on $E$ and consider the space $\A(\phi)$ of its extensions to states of $A$. This is a convex set, which is weakly $*$-closed in $\tilde \cS(A)$, so it is weakly $*$-compact as well. Hence $\A(\phi)$ has extreme points, which are necessarily extreme in $\tilde \cS(A)$. Indeed, if an extreme point $\psi \in A(\phi)$ can be written as $\psi = t \chi_1 + (1-t) \chi_2$ for some non-zero $\chi_1,\chi_2\in \tilde \cS(A)$ then also the restriction can be written as a convex combination: $\phi = t \chi_1|_E + (1-t) \chi_2|_E$. But since $\phi$ is pure,  this forces $\chi_1|_E = \chi_2|_E = \phi$ so that $\chi_1,\chi_2 \in \A(\phi)$, contradicting the fact that $\psi$ is extreme. 

Note that the proof also applies when replacing $E$ and $A$ by dense subspaces $\E$ and $\A$.
\endproof

We also record the following version of Jordan decomposition for matrix-ordered operator spaces, which is a direct consequence of the well-known $C^*$-algebraic version (see for instance \cite[Theorem 4.3.6]{KR83} for the case of unital operator systems). 
\begin{thm}[Jordan decomposition for matrix-ordered operator spaces]
  \label{thm:jordan}
   Let $E$ be a closed $*$-subspace in a $C^*$-algebra $A$ and  let $\E \subseteq \A$ be dense $*$-subspaces, so that $\E \subseteq E$ and $\A \subseteq A$. If there exists a norm-defining approximate order unit for $\A$ contained in $\E$,   
  then for each hermitian continuous linear functional $\phi :M_n(E) \to  \C$ on $E$ there exist positive linear functionals $\phi_+,\phi_- : M_n(E) \to \C$ such that $\phi = \phi_+ - \phi_-$ and $\| \phi\| = \| \phi_+ \| + \| \phi_-\|$. 
  \end{thm}
\proof
We extend $\phi:M_n(E) \to \C$ by the Hahn--Banach Theorem to a linear functional $\tilde \phi: M_n(A) \to \C$ so that $\tilde \phi|_{M_n(E)} = \phi$ and $\| \tilde \phi \| = \| \phi \|$. Here we may assume that $\tilde \phi$ is hermitian by replacing it by $\frac 12 (\tilde \phi + \tilde \phi^*)$ (which has the same norm as $\tilde \phi$ has). We then apply the Jordan decomposition ({\em cf.} \cite[Theorem 3.2.5]{Ped18}) for the $C^*$-algebra $M_n(A)$ to obtain positive functionals $\tilde \phi_\pm : M_n(A) \to \C$ with the property that
$$
\tilde \phi = \tilde \phi_+ - \tilde \phi_- ; \qquad \| \tilde \phi\| = \| \tilde \phi_+ \| + \| \tilde \phi_-\|.
$$
The restrictions $\phi_\pm := \tilde \phi_\pm |_{M_n(E)}$ are positive functionals on $M_n(E)$ and we have $\| \phi_\pm \| \leq \|\tilde \phi_\pm\|$. But since $\| \phi \| = \| \tilde \phi \|$ we find that for $\phi = \tilde \phi|_{M_n(E)} = \phi_+ - \phi_-$ we have
$$
\| \phi \| = \| \tilde \phi_+ \| + \| \tilde \phi_- \| \geq \| \phi_+ \| + \| \phi_- \|
$$
Finally, for $\phi = \tilde \phi|_{M_n(E)} = \phi_+ - \phi_-$ we also have $\| \phi \| = \| \phi_+ - \phi_- \| \leq  \| \phi_+ \| + \| \phi_- \|$, and this completes the proof.
\endproof

\begin{rem}
Note that uniqueness of the above Jordan decomposition of operator systems does not hold in general, unless $E = A$ is a unital $C^*$-algebra. See \cite[Theorem 4.3.6, Exercise 4.6.22]{KR83}.
  \end{rem}

In any case, the Jordan decomposition will allow us to say something about the Banach space dual of $M_n(E) \subseteq M_n(A)$ ---and consequently about the quasi-state space of $M_n(E)$--- under the assumption that there is an approximate order unit for $\A$ with the stated properties.  
Let us start with some notation.

As before we write $M_n(E)^*$ and $M_n(A)^*$ for the Banach space duals of $M_n(E)$ and $M_n(A)$, respectively. Also the notation $M_n(E)_h,M_n(A)_h$ will be used for the hermitian subspaces, as well as $M_n(E)_h^*, M_n(A)_h^*$ for the Banach space duals (where we note that $(M_n(E)^*)_h \cong (M_n(E)_h)^*$, isometrically isomorphic). There is a canonical order on the dual spaces which is respected by the restriction map by sending $M_n(A)^*_+ \to M_n(E)^*_+$.

We define as usual the annihilator of $M_n(E)$ as the closed subspace 
$$
M_n(E)^\perp := \{ \eta \in M_n(A)^* \mid \eta(x) = 0, \forall x\in M_n(E) \} \subseteq M_n(A)^*
$$
The quotient $M_n(A)^*/M_n(E)^\perp$ is a Banach space, on which we may introduce the following canonical order:
\begin{align}
  &(M_n(A)^*/M_n(E)^\perp)_+ := q( M_n(A)^*_+ ) \nonumber \\
  & \qquad = \{ \phi + M_n(E)^\perp \mid \phi \in M_n(A)^* \text{ and } \phi + \eta \geq 0 \text{ for some } \eta \in M_n(E)^\perp \}
  \label{eq:order-quot}
\end{align}
using the quotient map $q:  M_n(A)^* \to  M_n(A)^* /M_n(E)^\perp$.

\begin{rem}
The general construction of (matrix) orders on quotient spaces of operator systems is quite subtle (already in the unital case) and has been described in \cite{KPTT13}. 
\end{rem}

\begin{thm}
  \label{thm:quotient-dual}
   Let $E \subseteq A$ be a closed $*$-subspace in a $C^*$-algebra $A$ and let $\E \subseteq \A$ be dense $*$-subspaces with a norm-defining approximate order unit for $\A$ contained in $\E$.
  
Then the map $\rho:  M_n(A)^*_h /M_n(E)^\perp_h \to M_n(E)^*_h$ induced by restriction $\rho(\phi + M_n(E)^\perp_h) = \phi|_{M_n(E)_h}$ is an isometrical order isomorphism. 
\end{thm}
\proof
The map $\rho$ is clearly well-defined and injective. Moreover, for any $\eta \in M_n(E)^\perp_h$ we have
$$
\| \phi|_{M_n(E)_h} \| = \| (\phi + \eta)|_{M_n(E)_h} \| \leq \| \phi + \eta\|.
$$
Taking the infimum over all $\eta$ we find that $\rho$ is a contraction. Note also that $\rho$ is an order-preserving map. 

The inverse to $\rho$ is obtained by applying the Hahn--Banach theorem. Indeed, let a hermitian $\phi \in M_n(E)^*$ be given and write it in terms of a Jordan decomposition (Thm. \ref{thm:jordan}) as $\phi = \phi_+ - \phi_-$ where $\phi_\pm \in M_n(E)^*$ are two positive functionals such that $\| \phi \| = \| \phi_+\| + \| \phi_- \|$. By Proposition \ref{prop:pos-HB} there exist positive extensions $\tilde \phi_\pm \in M_n(A)^*$ with $\| \tilde \phi_\pm \| = \| \phi_\pm \|$. We set $\rho^{-1} (\phi) = \tilde \phi + M_n(E)^\perp$ where $\tilde \phi = \tilde \phi_+ - \tilde \phi_-$; this is a well-defined and linear map $M_n(E)^*_h \mapsto M_n(A)_h^*/M_n(E)_h^\perp$. For the norm we have
$$
\| \tilde \phi + M_n(E)^\perp_h \| \leq \| \tilde \phi \| \leq \|\tilde \phi_+\|+ \|\tilde \phi_- \|= \| \phi_+\|+ \|\phi_- \|= \|\phi \|,
$$ so that $\rho^{-1}$ is a contraction as well. Thus, both $\rho$ and $\rho^{-1}$ are isometries.

Finally, if $\phi= \phi_+$ is a positive functional then we have $\rho^{-1} (\phi) = \tilde \phi_+$ which is positive as well so that also $\rho^{-1}$ is order-preserving.
\endproof

\section{Operator systems, groupoid $C^*$-algebras, bonds and relations}
\label{sect:opsyst-bonds}
In this section we will introduce the main class of examples of operator systems which are associated to reflexive and symmetric relations on a set. It is convenient and also natural to formulate this in the slightly more general context of groupoids, and the corresponding $C^*$-algebras \cite{Ren80}.

We will consider a locally compact groupoid $G$ equipped with a (left invariant) Haar system $\nu = \{ \nu_x\}$. The space $C_c(G)$ of compactly supported complex-valued continuous functions on $G$ becomes a $*$-algebra with the convolution product and involution given by
$$
f \ast g (\gamma) = \int_{G_x} f(\gamma \gamma_1^{-1}) g(\gamma_1) d \nu_x(\gamma_1); \qquad f^*(\gamma) = \overline{f(\gamma^{-1})},
$$
where $x = s(\gamma)$ for any $\gamma \in G$. There is a norm given by
$$
\| f \|_1 = \sup_{x \in G^{(0)} } \left \{ \max ( \int_{G_x} |f(\gamma)| d\nu_x(\gamma) , \int_{G_x} |f(\gamma^{-1})| d\nu_x(\gamma) ) \right\} 
$$
and the completion $C_c(G)^{- \| \cdot \|_1}$ with respects to this norm turns out to be a Banach $*$-algebra. The universal $C^*$-algebra enveloping this Banach $*$-algebra will be called the {\em groupoid $C^*$-algebra}, and will be denoted by $C^*(G,\nu)$, or simply $C^*(G)$. We refer to \cite[Section II.1]{Ren80} for full details.  
\begin{rem}
  Note that for simplicity we have assumed that $G$ is Hausdorff, but this is not necessary. Indeed, arguing as in \cite{C82,KS02b} for the definition of $C^*(G)$ we should then replace $C_c(G)$ by a space of ---not necessarily continuous--- functions spanned by functions which vanish outside a compact Hausdorff subset $K \subset G$ and which are continuous on a neighborhood of $K$. 
\end{rem}

We now come to our main definition. 
\begin{defn}
  \label{defn:bonds}
  A {\em bond} is a 
  triple $(G,\nu,\Omega)$ consisting of a locally compact groupoid $G$, a Haar system $\nu = \{\nu_x\}$ and an open symmetric subset $\Omega \subseteq G$ containing the units $G^{(0)}$. 
  \end{defn}
When there is no ambiguity on the groupoid  and the Haar system,
we will also refer to the subset $\Omega$ as a bond. 

Then, if $\Omega_1,\Omega_2 \subseteq G$ are two bonds we have their composition, as usual:
$$
    \Omega_1 \circ \Omega_2 := \left \{ \gamma_1 \gamma_2 : \gamma_1 \in \Omega_1 ,\quad  \gamma_2 \in \Omega_2,\quad s(\gamma_1) = r(\gamma_2) \right\} 
$$
Clearly, this subset fails to be symmetric in general so we should consider instead the bond given by the symmetrized composition
$$
\Omega_1 * \Omega_2 := \Omega_1 \circ \Omega_2 \cup \Omega_2 \circ \Omega_1.
  $$

The following result follows directly from the definition of the involution in $C_c(G)$ and the assumption that a bond $\Omega$ is symmetric:
\begin{prop}
Let $(\Omega,G,\nu)$ be a bond. The closure of the subspace $C_c(\Omega) \subseteq C_c(G)$ in the $C^*$-algebra $C^*(G)$ is an operator system. 
  \end{prop}
We will denote this operator system as $E(\Omega,G,\nu)$ (or simply $E(\Omega)$ if no ambiguity can arise). Clearly, since $\supp(f\ast g) \subset \supp(f) \circ \supp(g)$ this subspace is in general not closed under the product in $C^*(G)$, unless the bond $\Omega$ is all of $G$ in which case we have $E(G)= C^*(G)$.

\begin{ex}
  Consider the infinite cyclic group $\Z$ as a groupoid. An example of a bond is given by an interval $\Omega_n = (-n,n) \subset \Z$ for some integer $n$. The corresponding operator system $E(\Omega_n, \Z)$ consists of sequences $a= (a_n)_{n \in \Z}$ in $C^*(\Z)$ with the following restricted support:
  $$
  a = ( \ldots 0, a_{-n+1}, a_{-n+2}, \ldots, a_{n-2}, a_{n-1} ,0 \ldots )
  $$
This operator system was analyzed at length in \cite{CS20}, where we called it the Fej\'er--Riesz operator system $C^*(\Z)_{(n)}$. In particular, it has propagation number $\infty$.
\end{ex}

\begin{ex}
  Consider now a {\em finite} cyclic group $C_m$. We can again consider the bond $\Omega_n = (-n,-n+1 \ldots, n-1,n) \subseteq \Z$ as in the previous example, but this time considered modulo $m$. We may then realize the corresponding operator system $E(\Omega_n, C_m)$ as the space of $m \times m$ banded circulant matrices with fixed band width (equal to $n$). In particular, one finds the propagation number to be finite.
\end{ex}

These two examples illustrate the important role played by the ambient groupoid, since for the same bond $(-n,n)$ but in different groups $\Z$ and $C_m$ one gets different, and not even Morita equivalent operator systems (in the sense of \cite{EKT21}). 

\subsection{Operator systems for tolerance relations}

An important class of examples of bonds is given by a symmetric and reflexive relation on a locally compact space $X$ equipped with a (Borel) measure of full support. Indeed, if we consider the pair groupoid $G = X \times X$ we may realize the relation as a symmetric subset $\cR \subseteq X \times X$ that contains the diagonal. In this case the operator system $E(\cR)$ is a subspace of the compact operators since $C^*(G) \cong \cK(L^2(X))$ \cite[Example 2.12]{Ren80}.

\begin{lma}
  \label{lma:prod-op-syst}
  Let $X$ be a locally compact space, equipped with a measure $\mu$ of full support. Then the operator system $E(\cR_1 \ast \cR_2)$ associated to the composition of two relations $\cR_1,\cR_2$ on $X$ coincides with $[E(\cR_1 )E(\cR_2) + E(\cR_2)E(\cR_1)]^{-\|\cdot\|}$.
\end{lma}
\proof 
Since $E(\cR_1 \ast \cR_2) = E(\cR_1 \circ \cR_2) + E(\cR_2 \circ \cR_1)$ it suffices to show that $E(\cR_1 \circ \cR_2) = [E(\cR_1)E(\cR_2)]^{-\|\cdot\|}$ as operator spaces.

Let us first show that $E(\cR_1)E(\cR_2) \subseteq E(\cR_1 \circ \cR_2)$. Indeed, for  kernels $F_j \in C_c(\cR_j)$ ($j=1,2$) the convolution product
$$
F_1 \ast F_2(x,y) = \int_X F_1(x,z) F_2(z,y) d\mu(z)
$$ is supported in $\cR_1 \circ \cR_2$, {\em i.e.} $F_1 \ast F_2 \in C_c(\cR_1 \circ \cR_2)$.  The result then follows by taking closures.

For the other inclusion take a function $F \in C_c(\cR_1\circ \cR_2 )$ with $\supp(F) =: K$. We may find compact subsets $K_j \subseteq \cR_j$ so that there exists an open subset $U$ with $K \subseteq U \subseteq K_1 \circ K_2$. We then take two non-negative functions $\chi_j \in C_c(\cR_j)$ with $\supp \chi_j = K_j$ and realize that their convolution product $\chi_1 \ast \chi_2$ will be a non-negative function in $C_c(\cR_1 \circ \cR_2)$ with support equal to $K_1 \circ K_2$ because $\mu$ has full support. It is thus strictly positive on $\supp( F)$ so that we may approximate $F$ by functions of the form 
$$
(x,y) \in \cR_1\circ \cR_2 \mapsto  \sum_{m=1}^M f_m (x) g_m(y) (\chi_1 \ast \chi_2)(x,y); \qquad (f_m,g_m \in C_c(X)).
$$
Moreover, kernels of this form satisfy
$$
G(x,y)= \sum_{m=1}^M f_m (x) g_m(y) (\chi_1 \ast \chi_2)(x,y)
=\sum_{m=1}^M G_{1,m} \ast G_{2,m} (x,y)
$$
where we have set $G_{1,m}(x,z) = f_m(x) \chi_1(x,z)$ and $G_{2,m}(z,y) = \chi_2(z,y)  g_m(y)$. Note that $G_{j,m}$ are continuous kernels with support in $K_j \subseteq \cR_j$ for all $m$ and we have thus approximated $F$ by functions in $C_c( \cR_1) \ast C_c(\cR_2)$. 
Taking closures then gives the desired inclusion $E(\cR_1 \circ \cR_2)\subseteq [E(\cR_1)E(\cR_2)]^{-\|\cdot\|}$.
\endproof

\subsection{Operator systems associated to tolerance relations on finite sets}
\label{sect:tol-finite}
We analyze in more detail the structure of the operator system $E(\cR)$ in the case that $\cR$ is a tolerance relation on a finite set $X$. In other words, we consider a bond $\cR$ on the pair groupoid associated to a finite set $X$. %
In this case, the unital operator system $E(\cR)$ is simply given by the linear span of rank-one operators $e_{xy}$ for all $(x,y) \in \cR$, acting linearly on the Hilbert space $\ell^2(X)$. We will denote by $\cK(\ell^2(X))$ the matrix $C^*$-algebra of all linear operators on $\ell^2(X)$. 

It is useful to view $\cR$ as a graph with the set of vertices given by $X$ and the set of edges by $\{ xy : (x,y) \in\cR\}$.

\begin{prop}
  \label{prop:tol-fin-prop}
Let $X$ be a finite set  and $\cR\subseteq X\times X$ a symmetric reflexive relation on $X$ and suppose that $\cR$ generates the full equivalence class $X \times X$ ({\em i.e.} the graph corresponding to $\cR$ is connected). Then 
  \begin{enumerate}
  \item The $C^*$-envelope of $E(\cR)$ is $C^*_\env(E(\cR)) = \cK(\ell^2(X))$.
    \item The propagation number can be expressed as $\text{prop}(E(\cR)) = \text{diam}(\cR)$, in terms of the graph diameter of $\cR$.
    \end{enumerate}
  \end{prop}
\proof
For any $x,y \in X$ with $d(x,y)=d$ there is a sequence $x=x_0 , x_1, \ldots, x_d =y$ of points in $X$ such that $(x_i,x_{i+1}) \in \cR$. Hence $e_{xy} = e_{x_0x_1} \cdots e_{x_{d-1}x_{d}}$ so that $C^*(E(\cR))= \cK(\ell^2(X))$. But then $C^*(E(\cR))$ is already the $C^*$-envelope because any  two-sided ideal is trivial and this applies in particular to the \v Silov ideal. We conclude that $\text{prop}(E(\cR)) \leq \text{diam}(\cR)$.

Conversely, suppose that $\text{prop}(E(\cR)) = n < \text{diam}(\cR)$ and take two points $x,y \in X$ with $d(x,y)>n$. Then for all $B_1, \ldots, B_n \in E(\cR)$ we have
$$
0 = \langle e_x, B_1 \cdots B_n e_y \rangle
$$
since $B_{i} (x_{i-1}, x_i) = 0$ whenever $d(x_{i-1},x_i) >1$ while $d(x_{i-1},x_i) \leq 1)$ would imply that $d(x,y)\leq n$. But then $\langle e_x, B e_y \rangle =0$ for all $B \in \cK(\ell^2(X))$ which is clearly a contradiction. We conclude that $\text{prop}(E(\cR)) = \text{diam}(\cR)$.
\endproof

For the description of the dual operator systems and state space, we make use of the following concrete realization of $E(\cR)$. 
First we identify $\cK(\ell^2(X)) \cong M_{|X|}(\C)$ using the canonical basis $\{e_x\}_{x\in X}$ of $\ell^2(X)$. Then we may write any element in $E(\cR)$ as a {\em sparse} $|X| \times |X|$ matrix with a fixed structure:
$$
E(\cR) \cong  \{ B \in \cK(\ell^2(X)) \mid B_{xy} = 0 \text{ if } (x,y) \notin \cR \}
\}.
$$
As a matter of fact, we may identify $E(\cR) \cong S_L(\cK(\ell^2(X)))$ 
using Schur--Hadamard (=entrywise) multiplication $S_L$ with the matrix associated to the relation $\cR$:
\begin{equation}
  \label{eq:schur}
L = (L_{xy}) ; \qquad L_{xy} = \left \{ \begin{array}{ll} 1 & \text{ if } (x,y) \in  \cR \\ 0 & \text{ if } (x,y)  \notin \cR  \end{array} \right.
\end{equation}
We will use the following basic result on binary matrices ({\em i.e.} matrices having entries in $\{0,1\}$).
\begin{lma}
  \label{lma:schur}
  Let $L$ be a binary matrix associated to a symmetric and reflexive relation. Then the following statements are equivalent:
  \begin{enumerate}
  \item[(i)]The matrix $L$ is positive semi-definite;
  \item[(ii)] The entries of the matrix $L$ satisfies the triangle inequalities;
    $$
L_{ij}+L_{jk}-L_{ik} \leq 1 \qquad (\forall i,j,k);
$$
\item[(iii)] The matrix $L$ is associated to an equivalence relation.
  \end{enumerate}
  \end{lma}
\proof
$(i) \implies (ii)$ follows since a violating triangle inequality should have $L_{ij}=L_{jk}=1$ and $L_{ik}=0$. This corresponds to a submatrix of the form
$$
\begin{pmatrix} 1 & 1 & 0 \\ 1 & 1 & 1 \\ 0 & 1 & 1 \end{pmatrix}
$$
which is indefinite, violating positive definiteness of $L$.

For $(ii) \implies (iii)$ one readily finds that the triangle inequalities imply transitivity for the relation represented by $L$.

Finally, for $(iii) \implies (i)$ we may  permute the columns and rows of the matrix $L$ to write it as a direct sum of block matrices containing only 1's. These blocks are all of the form $e^T e$ with $e = \begin{pmatrix} 1 & 1 & \cdots & 1 \end{pmatrix}$, and hence positive semi-definite.
\endproof

\begin{corl}
  \label{corl:schur}
  The matrix $L$ defined in Equation \eqref{eq:schur} is positive semi-definite if and only if the relation $\cR$ is an equivalence relation. In this case $E(\cR)$ is a matrix subalgebra of $\cK(\ell^2(X))$.
\end{corl}

We are now ready to apply the general results on dual operator systems from Section \ref{sect:jordan} to the finite-dimensional operator system $E(\cR) \subseteq \cK(\ell^2(X))$. We will make use of the dual pairing with $\cK(\ell^2(X))^d \cong \cK(\ell^2(X))$:
$$
\cK (\ell^2(X))^d \times \cK(\ell^2(X)) \to \C; \qquad (\rho, B) \mapsto \tr (\rho B)
$$

\begin{prop}
    \label{prop:state-space-finite}
  \begin{enumerate}
    \item 
  The dual operator system $E(\cR)^d$ of $E(\cR)$ is linearly isomorphic to $S_{L} (\cK(\ell^2(X)))$ with cones of positive elements given by
  $$
  M_n(E(\cR)^d)_+ = (S_{L})_n (M_n( \cK(\ell^2(X)))_+)
$$
where $S_{L}$ denote Schur multiplication with the matrix given in Equation \eqref{eq:schur}.
\item We have $S_L(M_n(E(\cR)^d)_+) =  S_L(M_n(E(\cR)^d))_+$ if and only if $\cR$ is an equivalence relation. 
  \end{enumerate}
  \end{prop}
\proof
As in Theorem \ref{thm:quotient-dual}  we have $E(\cR)^d \cong \cK(\ell^2(X))^d / E(\cR)^\perp $. Since 
$$
 E(\cR)^\perp \cong \{ \rho \in \cK(\ell^2(X)) \mid \rho_{xy} = 0 \text{ if } (x,y) \in \cR \} ,
 $$
 we find that the quotient map $\cK(\ell^2(X))^d \to \cK(\ell^2(X))^d /  E(\cR)^\perp $ is given by Schur multiplication with the matrix $L$. Hence  $E(\cR)^d \cong S_L(\cK(\ell^2(X)))$ as linear vector spaces.

 The order structure on $\cK(\ell^2(X))^d / E(\cR)^\perp $ is defined as in Equation \eqref{eq:order-quot} ({\em cf.} \cite{KPTT13}) in terms of the quotient map, so that $E(\cR)^d_+ \cong S_L(\cK(\ell^2(X))_+)$.  Since $\cK(\ell^2(X))$ is itself a matrix algebra, the extension to the matrix order is straightforward.

 The second claim follows from Lemma \ref{lma:schur}, in combination with the fact that $S_L$ is a completely positive map if and only if $L$ is positive semi-definite ({\em cf.}\cite[Theorem 3.7]{Pau02}). 
\endproof

For the pure states on the operator system $E(\cR)$ we have the following result. 
\begin{prop}
\label{prop:pure3}  Let $X$ be a finite set  and $\cR\subseteq X\times X$ a symmetric reflexive relation on $X$. The restriction of a pure state $\phi \in \mathbb P(\ell^2(X))$ to $E(\cR)$ is pure if and only if the restriction of the relation $\cR$ to the support $S$ of $\phi$ generates the full equivalence relation on $S$.
  \end{prop}
  \proof Let $v$ be a vector implementing the vector state $\phi$, one has $v_x\neq 0\iff x\in S$. Assume that the restriction of the relation $\cR$ to the support $S$ of $\phi$ generates the full equivalence relation on $S$.  Let then $w\in \ell^2(X)$ be a vector such that
  \begin{equation}  \label{inequ}
  \langle w, B w \rangle \leq \langle v, B v \rangle , \qquad \forall B \geq 0, B \in E(\cR)
 \end{equation}
  Applying \eqref{inequ} to diagonal matrices shows that $w_x=0$ for all $x\notin S$. Let then $(i,j)$, $i\neq j$, be a pair of elements of $S$ which belongs to the relation $\cR$. Consider the map $\iota:M_2(\mathbb C)\to E(\cR)$ which puts zeros outside the pair $(i,j)$  as a matrix with entries in $X$. The pure state on $M_2(\mathbb C)$ associated to $\phi\circ \iota$ is given by the vector $(v_i,v_j)$ up to normalization, it is pure because it is a vector state. By \eqref{inequ}  the inequality works for all positive elements of $M_2(\mathbb C)$ for the vector $(w_i,w_j)$. This implies that there exists a scalar $\lambda_{i,j}$ such that $(w_i,w_j)=\lambda_{i,j} (v_i,v_j)$. It follows that the ratio $w_i /v_i$ is independent of $i$ 
  since it is preserved by the relation $\cR$. \newline
  The converse follows from the existence of two disjoint equivalence classes in $S$ which split the restriction of the pure state $\phi$. \endproof

\section{Operator systems associated to metric spaces}
\label{sect:relations}

For us the motivating examples of tolerance relations ---and corresponding operator systems--- are given in terms of a metric space $(X,d)$ equipped with the relation
\begin{equation}
  \label{eq:rel-metric}
\cR_\eps := \left \{ (x,y) \in X \times X : d(x,y) < \epsilon \right\}
\end{equation}
Thus, we consider points in the space $X$ to be identical if they are within distance $\epsilon$. Clearly, this is reflexive and symmetric but not transitive (as long as $\eps$ is smaller than the diameter of $(X,d)$).

Let us analyze the composition rule for these relations, considering for instance $\cR_\eps \ast \cR_{\eps'}$ for different $\eps,\eps'>0$. Clearly, one expects this to be related to $\cR_{\eps+\eps'}$ but only for a certain class of metric spaces for which the metric distance between two points can be described in terms of smaller line segments. It turns out that the {\em path metric spaces} of \cite{Gro07} forms precisely the class of metric spaces that allows this (such spaces are called length spaces in \cite{Roe03}). 
\begin{defn}
  A {\em path metric space} is a metric space $(X,d)$ for which the distance between each pair of points equals the infimum of the lengths of curves joining the points.
\end{defn}
In particular, they satisfy the following two crucial, equivalent properties \cite[Theorem 1.8]{Gro07}:
\begin{enumerate}

\item for all $x,y \in X$ and $r>0$ there is a $z$ such that
  $$
\max \{ d(x,z), d(z,y) \} \leq \frac 12 d(x,y) + r
  $$
\item for all $x,y \in X$ and $r_1,r_2 \geq 0$ such that $r_1 + r_2 \leq d(x,y)$ one has
 \begin{equation}
   \label{eq:thm-gromov}
d\left( B_{r_1}(x), B_{r_2}(y) \right) \leq d(x,y)-r_1-r_2.
 \end{equation}
 \end{enumerate}
 We then arrive at the following result.
 \begin{lma}
   \label{lma:comp-rel}
   Let $(X,d)$ be a complete, locally compact path metric space and consider the relation $\cR_\eps$ defined in \eqref{eq:rel-metric}. For all $\eps,\eps'>0$ we have 
$$
  \cR_\eps * \cR_{\eps'}=  \cR_{\eps+\eps'}.
$$
\end{lma}
\proof
We clearly have $\cR_\eps \circ \cR_\eps' \subseteq \cR_{\eps+\eps'}$ which proves the inclusion $\cR_\eps * \cR_\eps' \subseteq \cR_{\eps+\eps'}$. For the converse, take $(x,y) \in \cR_{\eps+\eps'}$ so that $d(x,y) < \eps +\eps'$. Let $r_1+r_2 = d(x,y)$ be such that $r_1< \eps, r_2<\eps'$. Then with \eqref{eq:thm-gromov} we find that
$$
d\left( B_{r_1}(x), B_{r_2}(y) \right) =0.
$$
Hence there are sequences $(z_n)$ in $B_{r_1}(x)$ and $(z'_n)$ in $B_{r_2}(y)$ such that $d(z_n,z_n') \to 0$ as $n \to \infty$. By the Hopf--Rinow Theorem, the closed balls $\overline{B}_{r_1}(x)$ and $\overline{B}_{r_2}(y)$ are compact so that the sequences $(z_n), (z_n')$ have convergent sub-sequences that both converge to the same element $z \in X$. In fact, $z \in \overline{B}_{r_1}(x) \cap \overline{B}_{r_2}(y)$ so that $d(x,z) \leq r_1< \eps$ and $d(z,y) \leq r_2<\eps'$ as desired.
\endproof

\subsection{Operator systems for finite partial partitions}
In preparation for the study of the operator system $E(\cR_\eps)$ in the next subsection, we first analyze a class of finite-dimensional operator systems $E(\cR_\eps^P)$ which are associated to the relation $\cR_\eps$ and so-called finite partial partitions $P$ of a metric space $X$. 

\begin{defn}
  Let $X$ be a metric space and let $\eps>0$.
  \begin{itemize}
  \item A {\em finite partial $\eps$-partition} of $X$ is a finite collection $P = \{U_i\}$ of disjoint sets $U \subseteq X$ such that $\textup{diam}(U_i) < \eps$;
  \item We say that a finite partial partition $P'$ is a {\em refinement} of another finite partial partition $P$ if each set in $P'$ that lies in the support of $P$ is contained in a set in $P$ and if the supports $\cup_{U \in P} U \subseteq \cup_{U' \in P } U'$.
  \end{itemize}
  If $X$ is equipped with a (Borel) measure, then we will further assume that all sets in a finite partial $\eps$-partition are all measurable and have finite and non-zero measure.
  \end{defn}
Note that the notion of refinement induces a directed partial ordering on the set of finite partial $\eps$-partitions of a metric space. 

Let $X$ be a metric measure space and $P$ a finite partial $\eps$-partition. We can embed the finite-dimensional Hilbert space $\ell^2(P)$ in $L^2(X)$ by considering normalized characteristic function $1_U$ on $U$ as elements in $L^2(X)$. The corresponding finite-dimensional matrix algebra $\A_P$ is then defined by 
$$
\A_P =  \left\{ \sum_{U,V \in P} a_{UV} \ket{1_U} \bra{1_V} : a_{UV} \in \C \right\}
$$
and, in fact, we realize that $\A_P \cong \cK(\ell^2(P))$.
Note that the unit in $\A_P$ is given by $e_P = \sum_{U \in \cP}  \ket{1_U} \bra{1_U}$.
\begin{lma}
  \label{lma:struct-AP}
  Let $P,P'$ be finite partial $\eps$-partitions of a metric measure space $X$ and let $\A_P$ be as defined above.
  \begin{enumerate}
  \item If $P \leq P'$ then $\A_{P} \subseteq \A_{P'}$.
  \item The (algebraic) direct limit $\varinjlim \A_{P}$ is a dense $*$-subalgebra of the compact operators $\cK(L^2(X))$. Hence $\cK(L^2(X)) = \overline \varinjlim \A_{P}$ is the $C^*$-algebraic direct limit.
  \item The algebra $\varinjlim \A_{P}$ has a norm-defining approximate order unit (with respect to the matrix norm and matrix order induced from $\cK(L^2(X))$.
  \end{enumerate}
\end{lma}
\proof
(1) follows from the fact that each $U \in P$ can be written as a finite disjoint union of sets $U_j' \in P'$. Hence, $\ell^2(P) \subset \ell^2({P'})$ and thus also $\cK(\ell^2(P)) \subseteq \cK(\ell^2({P'}))$. 

For (2) we note that $\varinjlim \A_P$ is the linear span of rank-one operators of the form $\ket{1_U}\bra{1_V}$ with $\text{diam}(U),\textup{diam}(V)< \eps$. Now, any $L^2$-function on $X$ can be approximated by a continuous function, which in turn can be approximated by step functions of the form $\sum_{U \in P} c_U 1_U$ for some sufficiently fine $P$. As a consequence we find that $\varinjlim A_P$ is a dense subalgebra of the algebra $\cK_0(L^2(X))$ of all finite-rank operators, whose closure is $\cK(L^2(X))$.

For (3) we note that for any $x \in \A_{P}$ we have $- \| x \| e_P \leq x \leq \| x \| e_P $ and this indeed defines the norm of $x$. The approximate order unit for $\varinjlim \A_{P}$ is then given by $\{ e_P \}_P$, indexed by all finite partial $\eps$-partitions. 
\endproof

Next, we define relations $\cR_\eps^P$ 
associated to the partition $P$ and the relation $\cR_\eps$ introduced in Equation \eqref{eq:rel-metric}.

\begin{defn}
  Let $X$ be a metric measure space and $P$ a finite partial $\eps$-partition of $X$. We define a tolerance relation $\cR_\eps^P$ on the finite set $P$ by
  $$
\cR_\eps^P = \{ U \times V \mid \quad U,V \in P,  \quad U \times V \subseteq \cR_\eps \} \subseteq P \times P.
$$
\end{defn}
As in Section \ref{sect:tol-finite} the corresponding finite-dimensional operator system $E(\cR_\eps^P)$ is given by the linear span of operators of the form $\ket{1_U}\bra{1_V}$ with $U \times V \in  \cR_\eps$. 

\begin{lma}
  \label{lma:struct-EP}
  Let $P,P'$ be finite partial $\eps$-partitions of a 
  metric measure space $X$ and let $E(\cR_\eps^P)$ be as defined above.
  \begin{enumerate}    
 
  \item If $P \leq P'$ then $E(\cR_\eps^P) \subseteq E(\cR_\eps^{P'})$.
  \item The approximate order unit $\{ e_P \}_P$ of $\varinjlim A_P$ is contained in $E(\cR_\eps^P)$
  \end{enumerate}
\end{lma}
\proof
For (1) take $U,V \in P$ with $U \times V \in \cR_\eps$. As in the proof of Lemma \ref{lma:struct-AP} we may write $U = \cup_j U_j'$ and $V = \cup_k V_k'$ in terms of the refinement $P'$. But then also $U_j \times V_k \subseteq U \times V \subseteq \cR_\eps$ and hence the rank-one operator $\ket{1_U} \bra{1_V} \in E(\cR_\eps^P)$ can be written as $ \sum_{j,k}\ket{1_{U_j'}} \bra{1_{V_k'}}$ which is an element in $E(\cR_\eps^{P'})$. 

For (2) it suffices to note that $e_P = \sum_{U \in P} \ket{1_U}\bra{1_U}$ in terms of sets $U$ with diameter $< \eps$. But then $U \times U \subseteq \cR_\eps$ so that  $e_P \in E(\cR_\eps^P)$.
\endproof


As a special case of the results obtained at the end of Section \ref{sect:tol-finite} we have:

\begin{prop}
  \begin{enumerate}
  \item The operator system $E(\cR_\eps^P)$ is complete order-isomorphic to $S_L(\A_P) \subseteq \A_P $, in terms of Schur--Hadamard multiplication with the matrix
    $$
L = (L_{ij}) ; \qquad L_{ij} = \left \{ \begin{array}{ll} 1 & \text{ if } U_i \times U_j \subseteq \cR_\eps \\ 0 & \text{ if } U_i \times U_j \not \subseteq \cR_\eps  \end{array} \right.
    $$
\item 
    
  The dual operator system $E(\cR_\eps^P)^d$ of $E(\cR_\eps^P)$ is linearly isomorphic to $S_{L} (\A_P)$ with cones of positive elements given by
  $$
M_n(E(\cR_\eps^P)^d)_+ = (S_{L})_n (M_n(\A_P)_+) 
$$
\item We have $S_L((\A_P)_+) =  S_L(\A_P)_+$ if and only if $\cR_\eps^P$ is an equivalence relation.
  \end{enumerate}
\end{prop}
For the pure states of $E(\cR_\eps^P)$, Proposition \ref{prop:pure3} implies the following result.   
\begin{prop}
  The restriction of a pure state $\phi \in \P(l^2(P))$ to $E(\cR_\eps^P)$ is pure if and only if the restriction of $\cR_\eps^P$ to the support of $\phi$ generates the full equivalence relation.  
  \label{prop:pure-fin}
  \end{prop}

\begin{ex}
An interesting example of operator systems $E(\cR_\eps^P)$ are associated to finite $\eps$-partitions $P$ of the unit interval $[0,1)$. Let us consider partitions of the form
  \begin{equation}
    \label{eq:part-interval}
P = \left\{ U_k := \left[\frac { k} p , \frac{(k+1)}p \right) \right\}_{k=0}^{p-1}
\end{equation}
  We set a lower bound $\eps > 1/p$ so that each cell $U_k$ satisfies $U_k \times U_k \subseteq \cR_\eps$. The (unital) operator system $\E_P (\cR_\eps)$ can then be realized as the operator system $\E_{p,N}$ of $p \times p$ band matrices $B$ of band width $N$, {\em i.e.},
  $$
\E_{p,N} := \{  B = (b_{ij}) \in M_p(\C) \mid \qquad b_{ij} = 0 \text{ if } |i-j|>N \}
  $$
Indeed, the relation $\cR_\eps^P$ on $P$ that defines the operator system $E(\cR_\eps^P)$ can be written as
  \begin{equation}
    \label{eq:rel-RepsP}
  \cR_\eps^P = \{ U_k \times U_l : U_k  \times U_l \subset \cR_\eps \} =  \{ U_k \times U_l : |k-l| \leq N \} .
  \end{equation}
  so that a basis for $E(\cR_\eps^P)$ is thus given by $\{ \ket{1_{U_k}} \bra{1_{U_l}}\}_{|k-l| \leq N}$ which indeed spans the band matrices of the claimed band width.

  One quickly realizes that $\E_{p,N} \E_{p,N'} = \E_{p,N+N'}$ and that the propagation number of $\E_{p,N} \subseteq M_p(\C)$ is equal to $\lceil{p/N} \rceil$.

\end{ex}

\subsection{Operator systems for metric spaces at finite resolution}
Suppose now that $X$ is a path metric space, which is also a measure space with a measure of full support. 
Consider the operator systems $E(\cR_\eps) \subseteq \cK(L^2(X))$. It is convenient to consider it as a limit of the operator systems associated to finite partial partitions that we considered in the previous section. So let us write $\E(\cR_\eps) = \varinjlim E(\cR_\eps^P)$ for the algebraic direct limit considered in Lemma \ref{lma:struct-EP}; similarly we write $\A_\eps = \varinjlim \A_P$. 


\begin{prop}
  \label{prop:approx-unit}
  Let $X$ be a path metric measure space with a measure of full support. Let $\E(\cR_\eps) \subseteq \A_{\eps}$ be the (non-complete) operator system associated to the relation $\cR_\eps$ as defined above. Then
  \begin{enumerate}
  \item $\E(\cR_\eps)$ is a dense subspace of the operator system $E(\cR_\eps)$ and, consequently, we have $E(\cR_\eps) = \overline \varinjlim E(\cR_\eps^P)$, the closure of the algebraic direct limit;
  \item $\A_{\eps}$ is a dense $*$-subalgebra of the $C^*$-algebra $\cK(L^2(X))$;
  \item there exists a norm-defining approximate order unit for $\A_{\eps}$ which is contained in $\E(\cR_\eps)$.
    \end{enumerate}
  \end{prop}

\proof
For (1) we note that $E(\cR_\eps)$ is defined as the closure of the space of integral operators $\pi(F)$ on $L^2(X)$ with kernel $F \in C_c(\cR_\eps)$. 
We may cover the compact support of $F$ by a finite collection $\cP$ of open balls of radius $< \eps$. Hence, there are simple functions of the form $\sum_{U \times V \subseteq \cR_\eps} c_{UV} 1_{U \times V}$ that approximate $F$ in $L^\infty$-norm. But then the same is true in $L^2$-norm, again by compactness of the support of $F$. Consequently, the Hilbert--Schmidt operator $\pi(F)$ may be approximated by $\pi( 1_{U \times V}) = \ket{1_U}\bra{1_V} \in E(\cR_\eps^P)$.

The second and third statement were already established in Lemma \ref{lma:struct-EP}. 
\endproof

It is not surprising that the operator system $E(\cR_\eps)$ is closely related to the {\em Roe algebra} \cite{Roe93} associated to the coarse metric structure of $(X,d)$. Recall that in the present context of a metric measure space, the Roe algebra (or {\em translation $C^*$-algebra}) $C^*(X)$ is the norm closure of locally compact and finite propagation bounded operators on $L^2(X)$. Here an operator $T$ is called locally compact if $f T$ and $T f$ are compact for all $f \in C_0(X)$ and $T$ has propagation $<\eps$ means that for any $x,y \in X$ with $d(x,y)>\eps$ then $(x,y)$ is not in the support of $T$. One realizes as in \cite{OY15} that $C^*(X)$ is filtered by the collection of self-adjoint subspaces $(C^*(X)_\eps)_{\eps>0}$ where $C^*(X)_\eps$ is the closure of locally compact bounded operators with propagation less than $\eps$. Now, since the support of an operator $\pi(F) \in E(\cR_\eps)$ is given by the essential support of the kernel $F$, we conclude that $E(\cR_\eps) \subseteq C^*(X)_\eps$ for all $\eps>0$.

The following result shows the interplay between the propagation number for the operator system $E(\cR_\eps)$ and the finite propagation of the corresponding operators.

\begin{thm}
  \label{thm:prop-nr}
  Let $(X,d)$ be a complete, locally compact path metric measure space and $\mu$ a measure on $X$ with full support. Let  $E(\cR_\epsilon)$ be the operator system associated to the relation $d(x,y)< \epsilon$. Then  $$
  {\rm prop}(E(\cR_\epsilon)) = \lceil \delta/\epsilon\rceil
  $$
where $\lceil \cdot \rceil$ is the ceiling function and $\delta$ the diameter of $(X,d)$. Moreover the same equality holds after taking the minimal tensor product of $E(\cR_\epsilon)$ with the compact operators in a Hilbert space $H$.
\end{thm}
\proof

First observe that for $n \geq \lceil \delta/\epsilon\rceil$ we derive from Lemma \ref{lma:prod-op-syst} in combination with Lemma \ref{lma:comp-rel} that
\begin{equation}
  \label{eq:prod-op-syst-eps}
E(\cR_\eps)^{(n)} = E(\cR_{\eps} \ast \cdots * \cR_{\eps}) = E(\cR_{n \eps}) = \cK(L^2(X))
 \end{equation}
since $\cR_{n \eps} = X \times X$. Hence, 
we have $C^*({ E(\cR_\eps)}) \cong { \cK(L^2(X))}$. But then $C^*( {E(\cR_\eps)})$ is already the $C^*$-envelope of $E(\cR_\eps)$ because any boundary subsystem in the $C^*$-extension $C^*(E(\cR_\eps))$ of $E(\cR_\eps)$ is trivial and this applies in particular to the \v{S}ilov ideal. We conclude that ${\rm prop} (E(\cR_\eps)) \leq \lceil \delta/\epsilon \rceil$.

Now let $n$ be an integer such that $n\epsilon <\delta$ ({\em i.e.} $n < \lceil \delta/\eps \rceil$). Let $\xi,\eta\in L^2(X,\mu)$ be such that the distance between the essential supports of $\xi$ and $\eta$ is $>n\epsilon$. One has for $K_j = \pi(k_j) \in E(\cR_\eps)$ that 
$$
\langle \xi,K_1K_2\ldots K_n \eta \rangle=\int \overline{\xi(x_1)}k_1(x_1,x_2)\ldots k_n(x_n,x_{n+1})\eta(x_{n+1}) \prod d\mu(x_j)
$$
and the integrand vanishes identically since 
$$
d(x_j,x_{j+1})\leq \epsilon\  \forall j \Rightarrow 
d(x_1,x_{n+1})\leq n\epsilon \Rightarrow  \xi(x_1)\eta(x_{n+1})=0
$$ 
This shows that one has 
\begin{equation}
  \label{eq:inn-prod}
\langle \xi, T\eta\rangle =0, \qquad \forall T\in E(\cR_\epsilon)^{(n)}
\end{equation}
Thus $E(\cR_\epsilon)^{(n)}$ cannot be a $C^*$-algebra since by the argument leading to \eqref{eq:prod-op-syst-eps} it would then be the $C^*$-algebra of compact operators and the above vanishing does not hold for some compact operator. Hence ${\rm prop}(E(\cR_\eps)) \geq \lceil \delta/\eps \rceil$ which proves the first claim.

This argument can be extended to the minimal tensor product of $E(\cR_{\eps})$ with the compact operators $\cK(H)$ on a Hilbert space $H$ in the following way. Note first that for any integer $n>0$ we have $\left( E(\cR_\eps)\otimes_\min \cK(H)\right)^{(n)}  = E(\cR_\eps)^{(n)} \otimes_\min \cK(H)$ essentially because $\cK(H)^{(n)} = \cK(H)$ (in fact, $\cK(H)^{(n)}$ is a closed two-sided ideal in $\cK(H)$, hence exhausts all of $\cK(H)$). 

The next step consists of extending the inner product in \eqref{eq:inn-prod} to $E(\cR_\eps)^{(n)} \otimes_\min \cK(H)$ and apply the same argument as above. In fact, for any $\xi,\eta \in L^2(X,\mu)$ the map $T \to \langle \xi ,T \eta \rangle$ is completely bounded from $E(\cR_\epsilon)^{(n)}$ to $\C$, since it is the restriction of a linear functional on the $C^*$-algebra $\B(L^2(X))$. If we tensor it with the (also completely bounded) identity map on $\cK(H)$ we arrive at a completely bounded map
$$
E(\cR_\eps)^{(n)} \otimes_\min \cK(H) \to \cK(H) ; \qquad T \otimes A \mapsto \langle \xi, T\eta\rangle  \otimes A.
$$
We find that for $\xi,\eta \in L^2(X)$ with sufficiently separated essential supports as above, the $\cK(H)$-valued inner product is well-defined and in fact vanishes on all of $E(\cR_\eps)^{(n)} \otimes_\min \cK(H)$, showing that for $n < \lceil \delta/\eps \rceil$ we have $\left( E(\cR_\eps)\otimes_\min \cK(H)\right)^{(n)}  \neq \cK(L^2(X)) \otimes_\min \cK(H)$.
\endproof

This is in line with the result obtained for tolerance relations on finite sets in Proposition \ref{prop:tol-fin-prop}. In fact, that result can be derived as a special case from Theorem \ref{thm:prop-nr}. In order to see this, equip the finite set $X$ with the graph metric obtained from the given relation $\cR \subseteq X \times X$, then for any $1 < \eps < 2$ one readily sees that $\cR_\eps  = \cR$. 

\bigskip

On physical grounds one could expect a relation between the operator systems $E(\cR_\eps^{S^1})$ for the circle at finite resolution $\eps$ and the Toeplitz operator system $\Toep{n}$ obtained by spectral projection of the circle onto the first $n$ Fourier modes. The latter operator system has been considered at length in our previous work \cite{CS20} and also in \cite{Far21}. The reason for having such a relation is that the truncation in momentum space is expected to correspond to introducing a finite resolution in position space. However, in view of the above results we conclude that such a relation cannot be found directly by means of a complete order-isomorphism between the operator systems, and in fact even a stronger result is true.

\begin{corl} For $\epsilon$ small enough, there is no value of $n$ such that the operator systems
$\Toep{n}$ and  $E(\cR^{S^1}_\eps)$ are stably equivalent.
\end{corl}
\proof
If was shown in \cite[Proposition 4.2]{CS20} that ${\rm prop}(\Toep{n})=2$ for any $n$, while the fact that the propagation number is invariant up to stable equivalence is \cite[Proposition 2.42]{CS20}.
\endproof

\subsubsection{State space of $E(\cR_\eps)$}
\label{sect:states}

It is well-known that $\cK^* \cong \L^1$ isometrically via the pairing $(\rho, T) \in \L^1\times \cK \mapsto \tr \rho T$. This duality respects the natural ordering by positive elements on both sides in the sense that
\begin{align*}
  \rho \in \L^1_+ \iff \tr \rho T \geq 0 , \qquad \forall T \in \cK_+,\\
  \intertext{and also}
  T \in \cK_+ \iff \tr \rho T \geq 0 , \qquad \forall \rho \in \L^1_+.
\end{align*}
We then have as a special case of Theorem \ref{thm:quotient-dual}.
\begin{prop}
\label{prop:isom-quot-cpt}
  There is an isometrical order isomorphism $E(\cR_\eps)^*_h \cong \L^1_h/ E (\cR_\eps)^\perp_h$. 
  \end{prop}

\begin{corl}
Any state on $E(\cR_\eps) \subseteq \cK(L^2(X))$ can be written as $T \mapsto \tr \rho T$ ($T \in E(\cR_\eps)$), in terms of a trace-class density operator $\rho = \sum_i \lambda_i \ket{\xi_i}\bra{\xi_i}$ with $\lambda_i \geq 0$ and $\xi_i \in L^2(X)$.
  \end{corl}

We see from Proposition \ref{prop:ext-pure} that for any pure state $\phi$ on $E(\cR_\eps)$, there exists a pure state on $\cK(L^2(X))$ that extends $\phi$. Since pure states on the compacts are given by vector states, we find that any pure state $\varphi$ on $E(\cR_\eps)$ can be written as
$$
\varphi: A \mapsto \langle \psi,A \psi \rangle_{L^2(X)}
$$
for some $\psi \in L^2(X)$ of norm one. In general, the correspondence between $\mathcal P(E(\cR_\eps))$ and $\mathbb P(L^2(X))$ is not expected to be one to one. However, this happens if one restricts to pure states on $\cK(L^2(X))$ whose support is $\eps$-connected in the following sense.
\begin{defn}
  A subset $Y$ of a metric space $(X,d)$ is $\eps$-connected if and only if the equivalence relation on $Y$ generated by the relation $d(y,y')<\eps$ contains only one equivalence class.
\end{defn}
We let $\mathbb P_\eps(L^2(X))\subset\mathbb P(L^2(X))$ be the subset of $\mathbb P(L^2(X))$ formed of pure states $\vert \psi\rangle\langle\psi\vert$ where the essential support of $\psi$ is $\eps$-connected. 
   \begin{lma} \label{lma:pure1} The restriction map $\rho:\mathbb P_\eps(L^2(X))\to \cS(E(\cR_\eps))$ from $\eps$-connected pure states  on $\cK(L^2(X))$ to states on $E(\cR_\eps)$ is injective.   	
   \end{lma}
   \proof The state $\rho(\vert \psi\rangle\langle\psi\vert)$ is given by the pairing
   $$
   \rho(\vert \psi\rangle\langle\psi\vert)(k)=\int k(x,y)\bar\psi(x)\psi(y)d\mu(x)d\mu(y)
   $$
   which determines the product $\bar\psi(x)\psi(y)$ almost everywhere, for $d(x,y)<\eps$. This shows that if $\rho(\vert \psi\rangle\langle\psi\vert)=\rho(\vert \psi'\rangle\langle\psi'\vert)$, then one has $\vert \psi(x)\vert=\vert \psi'(x)\vert$ almost everywhere and there exists a measurable map $u:X\to U(1)$ well determined on the  support of $\psi$, such that $\psi'(x)=u(x)\psi(x)$. Let us show that $u$ is constant. One has
   $$
   d(x,y)<\eps\Rightarrow \bar\psi(x)\psi(y)=\bar\psi'(x)\psi'(y)
   $$
and this shows that for $x,y$ in the essential support of $\psi$ one has 
$$
   d(x,y)<\eps\Rightarrow u(x)=u(y).
   $$
Thus since the support of $\psi$ is $\eps$-connected, one gets that $u$ is constant. \endproof 
\begin{lma}
  \label{lma:pure2}
  Let $X$ be a metric space with a measure $\mu$ with full support.   
  The restriction of a pure state  $\psi\in\mathbb P(L^2(X))$ to $E(\cR_\eps)$ is pure if and only if $\psi$ is $\eps$-connected.   	
   \end{lma}
   \proof Assume first that $\psi$ is not $\eps$-connected. Let $S$ be the essential support of $\psi$ and $S=S_1\cup S_2$ a partition of $S$ in two subsets such that the distance $d(S_1,S_2)\geq \epsilon$. With $\psi=\vert \xi\rangle\langle\xi\vert$, let $\xi_j=1_{S_j}\xi$ so that $\xi=\xi_1+\xi_2$ and one has 
   $$
   d(S_1,S_2)\geq \epsilon \Rightarrow \langle T\xi_i,\xi_j\rangle=0, \  \forall i\neq j, \  \forall T\in E(\cR_\eps).
   $$
   One thus gets that 
   $$
   \psi(T)= \langle T\xi_1,\xi_1\rangle+ \langle T\xi_2,\xi_2\rangle, \  \forall T\in E(\cR_\eps),
   $$
   and one gets that the restriction of $\psi$ to $E(\cR_\eps)$ is not pure since it is the non-trivial convex combination of the normalized states $\vert \xi_j\rangle\langle\xi_j\vert\times \Vert \xi_j\Vert^{-2}$. \newline

   For the converse, assume that the essential support $S$ of $\xi$ is $\eps$-connected. The vector state defined by $\xi$ is pure if we have
     $$
\langle \eta, T \eta \rangle \leq \langle \xi, T \xi \rangle ; \qquad \forall T \geq 0, T \in E(\cR_\eps) \implies \eta \in \C \xi.
$$
Consider then such an $\eta$ and note that the inequality holds for any positive $T \in E(\cR_\eps^P) \subset E(\cR_\eps)$ for any finite partial $\eps$-partition $P$.

For any measurable subset $U \subset X \setminus S$ of sufficiently small diameter but with non-zero measure, we may consider the finite partial $\eps$-partition $P = \{ U \}$. Since we have $e_P \in E(\cR_\eps^P)$ we find that $\langle \eta, e_P \eta \rangle \leq \langle \xi, e_P \xi \rangle =0$ . In other words, $\int_U |\eta|^2 \leq \int_U | \xi|^2 =0 $ and we may conclude that the essential support of $\eta$ is contained in $S$.

For any two $x,y \in S$ so that $d(x,y) <\eps$ we may consider sets $U,V$ of diameter $< \eps$ with $x \in U, y \in V$ and $U \times V \subset \cR_\eps$. By the assumption on the support of $\mu$ they may be taken to have non-zero measure. Then we find for the partial $\eps$-partitions $P = \{ U,V\}$ that the vector state restricted to $E(\cR_\eps^P)$ becomes $e_P \xi =:( \xi_U , \xi_V) \in \C^2$. By Proposition \ref{prop:pure-fin} this is a pure state on $E(\cR_\eps^P)$, whence the vector $e_P \eta$ is a complex multiple of $e_P \xi$. 

   Let us now fix $V$ and suppose that $\eta_V = c \xi_V$ for some $c$ where we assume that $\xi_V \neq 0$. Now for any $U' \subset U$ we still have $U' \times V \subset \cR_\eps$ so that the argument from the previous paragraph applies and we find that $\eta_{U'} = c \xi_{U'}$. We conclude that $\eta = c \eta$ in a neighborhood of $x$.  

   Repeating this argument to any point at distance $< \eps$ from $x$ and using $\eps$-connectedness of $S$ we find that $\eta$ is a constant multiple of $\xi$ on the essential support of $\xi$. 
\endproof
  
\begin{thm}
\label{prop:pure2} Let $X$ be a metric space with a measure $\mu$ with full support.   The restriction map $\rho:\mathbb P_\eps(L^2(X))\to \mathbb \cP(E(\cR_\eps))$ from $\eps$-connected pure states  on $\cK(L^2(X))$ to pure states on $E(\cR_\eps)$ is bijective. 
  \end{thm}
\proof
First Lemma \ref{lma:pure2} shows that $\rho$ maps $\eps$-connected pure states to pure states. Lemma \ref{lma:pure1} shows that $\rho$ is injective. Let us show that it is surjective. By Hahn--Banach extension there exists, given a pure state $\phi$ on $ E(\cR_\eps)$ a pure state $\psi$ on $\cK(L^2(X))$ whose restriction gives $\phi$. By Lemma \ref{lma:pure2} the pure state $\psi$ is $\eps$-connected. \endproof

\newcommand{\noopsort}[1]{}\def\cprime{$'$}

\end{document}